\documentclass[11pt, reqno]{amsart}
\usepackage{amssymb}
\usepackage{amsmath}
\usepackage{amsthm}
\usepackage{amsfonts}
\usepackage[latin5]{inputenc}

\input{tcilatex}

\begin{document}
\title[Fourth Fundamental Form and  $i$-th Curvature Formulas in ${\mathbb{E}%
}^{4}$]{Fourth Fundamental Form and \\
$i$-th Curvature Formulas in ${\mathbb{E}}^{4}$}
\author{Erhan Güler}
\address{Erhan Güler: Bart\i n University, Faculty of Sciences, Department
of Mathematics, 74100 Bart\i n, Turkey.}
\email{eguler@bartin.edu.tr }
\urladdr{https://orcid.org/0000-0003-3264-6239}
\date{Received: October 30, 2020}
\thanks{This paper is in final form and no version of it will be submitted
for publication elsewhere.}

\begin{abstract}
We introduce fourth fundamental form $IV,$ and $i$-th curvature formulas of
hypersurfaces in the four dimensional Euclidean geometry ${\mathbb{E}}^{4}$.
Defining fourth fundamental form and $i$-th curvatures for hypersurfaces, we
calculate them on rotational hypersurface. In addition we study rotational
hypersurface satisfying $\Delta ^{IV}\mathbf{x=Ax}$ for some $4\times 4$
matrix $\mathbf{A.}$
\end{abstract}

\maketitle

\section{Introduction}

Refering Chen \cite{C0, C, C2, C3}, the researches of submanifolds of finite
type whose immersion into $\mathbb{E}^{m}$ (or $\mathbb{E}_{\nu }^{m}$) by
using a finite number of eigenfunctions of their Laplacian has been studied
by geometers for almost a half century.

Takahashi \cite{T} proved that a connected Euclidean submanifold is of
1-type, iff it is either minimal in $\mathbb{E}^{m}$ or minimal in some
hypersphere of $\mathbb{E}^{m}$. Submanifolds of finite type closest in
simplicity to the minimal ones are the 2-type spherical submanifolds (where
spherical means into a sphere). Some results of 2-type spherical closed
submanifolds were obtained by \cite{BarrosChen, BarrosGaray, C}. Garay gave 
\cite{Garay2} an extension of Takahashi's theorem in $\mathbb{E}^{m}$. Cheng
and Yau worked hypersurfaces with constant scalar curvature; Chen and
Piccinni \cite{CP} studied submanifolds with finite type Gauss map in $%
\mathbb{E}^{m}$. Dursun \cite{Dursun} gave hypersurfaces with pointwise
1-type Gauss map in $\mathbb{E}^{n+1}$.

Considering ${\mathbb{E}}^{3}$; Takahashi \cite{T} stated that minimal
surfaces and spheres are the only surfaces satisfying the condition $\Delta
r=\lambda r,$ $\lambda \in \mathbb{R}$; Ferrandez, Garay and Lucas \cite{FGL}
proved that the surfaces satisfying $\Delta H=AH$, $A\in Mat(3,3)$ are
either minimal, or an open piece of sphere or of a right circular cylinder;
Choi and Kim \cite{CK} characterized the minimal helicoid in terms of
pointwise 1-type Gauss map of the first kind; Garay \cite{Garay} worked a
certain class of finite type surfaces of revolution; Dillen, Pas and
Verstraelen \cite{DPV} proved that the only surfaces satisfying $\Delta
r=Ar+B,$ $A\in Mat(3,3),$ $B\in Mat(3,1)$ are the minimal surfaces, the
spheres and the circular cylinders; Stamatakis and Zoubi \cite%
{StamatakisZoubi} considered surfaces of revolution satisfying $\Delta
^{III}x=Ax$; Senoussi and Bekkar \cite{SB} studied helicoidal surfaces $%
M^{2} $ which are of finite type with respect to the fundamental forms $I,II$
and $III,$ i.e., their position vector field $r(u,v)$ satisfies the
condition $\Delta ^{J}r=Ar,$ $J=I,II,III,$ where $A\in Mat(3,3)$; Kim, Kim
and Kim \cite{Kim et al} focused Cheng-Yau operator and Gauss map of
surfaces of revolution.

General rotational surfaces in 4-space were originated by Moore \cite{M, M2}%
. \ Focusing on ${\mathbb{E}}^{4}$; Hasanis and Vlachos \cite{HV} considered
hypersurfaces with harmonic mean curvature vector field; Cheng and Wan \cite%
{CW} considered complete hypersurfaces with $CMC$; Kim and Turgay \cite%
{KimTurgay} studied surfaces with $L_{1}$-pointwise $1$-type Gauss map;
Arslan et. al. \cite{Arslanetal} worked Vranceanu surface with pointwise $1$%
-type Gauss map; Kahraman Aksoyak and Yayl\i\ \cite{AksoyakYayli} introduced
flat rotational surfaces with pointwise 1-type Gauss map; Güler, Magid and
Yayl\i\ \cite{GulerMagidYayli} considered the helicoidal hypersurfaces; Gü%
ler, Hac\i salihoğlu and Kim \cite{GulerHacisalihogluKim} studied Gauss map
and the third Laplace-Beltrami operator of the rotational hypersurface; Gü%
ler and Turgay \cite{GulerTurgay} introduced Cheng-Yau operator and Gauss
map of rotational hypersurfaces; Güler \cite{Guler1} studied rotational
hypersurfaces satisfying $\Delta ^{I}R=AR$, where $A\in Mat(4,4).$

In Minkowski 4-space $\mathbb{E}_{1}^{4}$; Ganchev and Milousheva \cite{GM}
considered analogue of surfaces of \cite{M, M2}; Arvanitoyeorgos,
Kaimakamais and Magid \cite{AKM} showed that if the mean curvature vector
field of $M_{1}^{3}$ satisfies the equation $\Delta H=\alpha H$ ($\alpha $ a
constant), then $M_{1}^{3}$ has $CMC$; Arslan and Milousheva worked meridian
surfaces of elliptic or hyperbolic type with pointwise 1-type Gauss map;
Turgay gave some classifications of Lorentzian surfaces with finite type
Gauss map; Dursun and Turgay studied space-like surfaces in with pointwise
1-type Gauss map. Kahraman Aksoyak and Yayl\i\ \cite{AksoyakYayli2} focused
general rotational surfaces with pointwise 1-type Gauss map in $\mathbb{E}%
_{2}^{4}$. Bektaş, Canfes and Dursun \cite{BCD} classified surfaces in a
pseudo-sphere with 2-type pseudo-spherical Gauss map in $\mathbb{E}_{2}^{5}.$

In literature, there is no any work about fourth fundamental form $f_{ij}$
(i.e. $IV$\textbf{) }and $i$-th curvature formulas $\mathfrak{C}_{i}$, where 
$i=0,\hdots,3,$ of rotational hypersurface in the four dimensional Euclidean
space $\mathbb{E}^{4}$.

We consider fourth fundamental form $IV,$ and $i$-th curvature formulas $%
\mathfrak{C}_{i}$ of hypersurfaces in the four dimensional Euclidean
geometry ${\mathbb{E}}^{4}$. In Section 2, we give some basic notions of the
four dimensional Euclidean geometry. Defining fourth fundamental form and $i$%
-th curvature for hypersurfaces, we calculate $\mathfrak{C}_{i}$ and fourth
fundamental form of rotational hypersurface in Section 3. Finally, in the
last section, we study rotational hypersurface satisfying $\Delta ^{IV}%
\mathbf{x=Ax}$ for some $4\times 4$ matrix $\mathbf{A}$ in $\mathbb{E}^{4}$.

\section{Preliminaries}

In this section, giving some of basic facts and definitions, we describe
notations used whole paper. Let $\mathbb{E}^{m}$ denote the Euclidean $m$%
-space with the canonical Euclidean metric tensor given by $\widetilde{g}%
=\langle \ ,\ \rangle =\sum\limits_{i=1}^{m}dx_{i}^{2},$ where $(x_{1},x_{2},%
\hdots,x_{m})$ is a rectangular coordinate system in $\mathbb{E}^{m}$.
Consider an $m$-dimensional Riemannian submanifold of the space $\mathbb{E}%
^{m}$. We denote the Levi-Civita connections of $\mathbb{E}^{m}$ and $M$ by $%
\widetilde{\nabla }$ and $\nabla $, respectively. We shall use letters $%
X,Y,Z,W$ (resp., $\xi ,\eta $) to denote vectors fields tangent (resp.,
normal) to $M$. The Gauss and Weingarten formulas are given, respectively,
by 
\begin{eqnarray}
\widetilde{\nabla }_{X}Y &=&\nabla _{X}Y+h(X,Y),  \label{MEtomWeingarten} \\
\widetilde{\nabla }_{X}\xi &=&-A_{\xi }(X)+D_{X}\xi ,
\end{eqnarray}%
where $h$, $D$ and $A$ are the second fundamental form, the normal
connection and the shape operator of $M$, respectively.

For each $\xi \in T_{p}^{\bot }M$, the shape operator $A_{\xi }$ is a
symmetric endomorphism of the tangent space $T_{p}M$ at $p\in M$. The shape
operator and the second fundamental form are related by 
\begin{equation*}
\left\langle h(X,Y),\xi \right\rangle =\left\langle A_{\xi }X,Y\right\rangle
.
\end{equation*}%
The Gauss and Codazzi equations are given, respectively, by 
\begin{eqnarray}
\langle R(X,Y,)Z,W\rangle &=&\langle h(Y,Z),h(X,W)\rangle -\langle
h(X,Z),h(Y,W)\rangle ,  \label{MinkCodazzi} \\
(\bar{\nabla}_{X}h)(Y,Z) &=&(\bar{\nabla}_{Y}h)(X,Z),
\end{eqnarray}%
where $R,\;R^{D}$ are the curvature tensors associated with connections $%
\nabla $ and $D$, respectively, and $\bar{\nabla}h$ is defined by 
\begin{equation*}
(\bar{\nabla}_{X}h)(Y,Z)=D_{X}h(Y,Z)-h(\nabla _{X}Y,Z)-h(Y,\nabla _{X}Z).
\end{equation*}

\subsection{Hypersurfaces of Euclidean space}

Now, let $M$ be an oriented hypersurface in the Euclidean space $\mathbb{E}%
^{n+1}$, $\mathbf{S}$ its shape operator (i.e. Weingarten map) and $x$ its
position vector. We consider a local orthonormal frame field $\{e_{1},e_{2},%
\hdots,e_{n}\}$ of consisting of principal directions of $M$ corresponding
from the principal curvature $k_{i}$ for $i=1,2,\hdots n$. Let the dual
basis of this frame field be $\{\theta _{1},\theta _{2},\hdots,\theta _{n}\}$%
. Then the first structural equation of Cartan is 
\begin{equation}
d\theta _{i}=\sum\limits_{i=1}^{n}\theta _{j}\wedge \omega _{ij},\quad
i,j=1,2,\hdots,n,  \label{CartansFirstStructural}
\end{equation}%
where $\omega _{ij}$ denotes the connection forms corresponding to the
chosen frame field. We denote the Levi-Civita connection of $M$ and $\mathbb{%
\ E}^{n+1}$ by $\nabla $ and $\widetilde{\nabla }$, respectively. Then, from
the Codazzi equation \eqref{MinkCodazzi}, we have 
\begin{eqnarray}
e_{i}(k_{j}) &=&\omega _{ij}(e_{j})(k_{i}-k_{j}),  \label{Coda2} \\
\omega _{ij}(e_{l})(k_{i}-k_{j}) &=&\omega _{il}(e_{j})(k_{i}-k_{l})
\end{eqnarray}%
for distinct $i,j,l=1,2,\hdots,n$.

We put $s_{j}=\sigma _{j}(k_{1},k_{2},\hdots,k_{n})$, where $\sigma _{j}$ is
the $j$-th elementary symmetric function given by 
\begin{equation*}
\sigma _{j}(a_{1},a_{2},\hdots,a_{n})=\sum\limits_{1\leq i_{1}<i_{2}<\hdots%
,i_{j}\leq n}a_{i_{1}}a_{i_{2}}\hdots a_{i_{j}}.
\end{equation*}%
We use following notation 
\begin{equation*}
r_{i}^{j}=\sigma _{j}(k_{1},k_{2},\hdots,k_{i-1},k_{i+1},k_{i+2},\hdots%
,k_{n}).
\end{equation*}%
By the definition, we have $r_{i}^{0}=1$ and $s_{n+1}=s_{n+2}=\cdots =0$. We
call the function $s_{k}$ as the $k$-th mean curvature of $M$. We would like
to note that functions $H=\frac{1}{n}s_{1}$ and $K=s_{n}$ are called the
mean curvature and Gauss-Kronecker curvature of $M$, respectively. In
particular, $M$ is said to be $j$-minimal if $s_{j}\equiv 0$ on $M$.

In $\mathbb{E}^{n+1},$ to find the $i$-th curvature formulas $\mathfrak{C}%
_{i}$ (Curvature formulas sometimes are represented as mean curvature $%
H_{i}, $ and sometimes as Gaussian curvature $K_{i}$ by different writers,
such as \cite{AliasGurbuzGeomDed2006} and \cite{K}. We will call it just $i$%
-th curvature $\mathfrak{C}_{i}$ in this paper.), where $i=0,..,n,$ firstly,
we use the characteristic polynomial of $\mathbf{S}$\textbf{:} 
\begin{equation}
P_{\mathbf{S}}(\lambda )=0=\det (\mathbf{S}-\lambda
I_{n})=\sum\limits_{k=0}^{n}\left( -1\right) ^{k}s_{k}\lambda ^{n-k},
\label{det}
\end{equation}%
where $i=0,..,n,$ $I_{n}$ denotes the identity matrix of order $n.$ Then, we
get curvature formulas $\binom{n}{i}\mathfrak{C}_{i}=s_{i},$. That is, $%
\binom{n}{0}\mathfrak{C}_{0}=s_{0}=1$ (by definition), $\binom{n}{1}%
\mathfrak{C}_{1}=s_{1},\hdots,\binom{n}{n}\mathfrak{C}_{n}=s_{n}=K.$

$k$-th fundamental form of $M$ is defined by $I\left( \mathbf{S}^{k-1}\left(
X\right) ,Y\right) =\left\langle \mathbf{S}^{k-1}\left( X\right)
,Y\right\rangle .$ So, we have%
\begin{equation}
\sum\limits_{i=0}^{n}\left( -1\right) ^{i}\binom{n}{i}\mathfrak{C}%
_{i}I\left( \mathbf{S}^{n-i}\left( X\right) ,Y\right) =0.
\label{Fundamentalforms}
\end{equation}%
In particular, one can get classical result $\mathfrak{C}_{0}III-2\mathfrak{C%
}_{1}II+\mathfrak{C}_{2}I=0$ of surface theory for $n=2.$ See \cite{K} for
details.

For a Euclidean submanifold $x$: $M\longrightarrow \mathbb{E}^{m}$, the
immersion $\left( M,x\right) $ is called \textit{finite type,} if $x$ can be
expressed as a finite sum of eigenfunctions of the Laplacian $\Delta $ of $%
\left( M,x\right) ,$ i.e. $x$ $=$ $x_{0}+\sum_{i=1}^{k}x_{i}$, where $x_{0}$
is a constant map, $x_{1},\hdots,x_{k}$ non-constant maps, and $\Delta
x_{i}=\lambda _{i}x_{i},$ $\lambda _{i}\in \mathbb{R}$, $i=1,\hdots,k$. If $%
\lambda _{i}$ are different, $M$ is called \textit{k-type}. See \cite{C} for
details.

\subsection{Rotational hypersurfaces}

We will obtain a rotational hypersurface (rot-hypface for short) in
Euclidean 4-space. Before we proceed, we would like to note that the
definition of rot-hypfaces in Riemannian space forms were defined in \cite%
{CDRotHyp}. A rot-hypface $M\subset \mathbb{E}^{n+1}$ generated by a curve $%
\mathcal{C}$ around an axis $\mathcal{C}$ that does not meet $\mathcal{C}$
is obtained by taking the orbit of $\mathcal{C}$ under those orthogonal
transformations of $\mathbb{E}^{n+1}$ that leaves $\mathfrak{r}$ pointwise
fixed (See \cite[Remark 2.3]{CDRotHyp}).

Throughout the paper, we shall identify a vector ($a,b,c,d$) with its
transpose. Consider the case $n=3$, and let $\mathcal{C}$ be the curve
parametrized by%
\begin{equation}
\gamma (u)=\left( f(u),0,0,\varphi \left( u\right) \right) .  \label{gamma}
\end{equation}%
If $\mathfrak{r}$ is the $x_{4}$-axis, then an orthogonal transformations of 
$\mathbb{E}^{n+1}$ that leaves $\mathfrak{r}$ pointwise fixed has the form 
\begin{equation*}
\mathbf{Z}(v,w)=\left( 
\begin{array}{cccc}
\cos v\cos w & -\sin v & -\cos v\sin w & 0 \\ 
\sin v\cos w & \cos v & -\sin v\sin w & 0 \\ 
\sin w & 0 & \cos w & 0 \\ 
0 & 0 & 0 & 1%
\end{array}%
\right) ,\text{ \ }v,w\in {\mathbb{R}}.
\end{equation*}%
Therefore, the parametrization of the rot-hypface generated by a curve $%
\mathcal{C}$ around an axis $\mathfrak{r}$ is given by 
\begin{equation}
\mathbf{x}(u,v,w)=\mathbf{Z}(v,w)\gamma (u).  \label{1}
\end{equation}

Let $\mathbf{x}=\mathbf{x}(u,v,w)$ be an isometric immersion from $%
M^{3}\subset \mathbb{E}^{3}$ to $\mathbb{E}^{4}$. Triple vector product of $%
\overrightarrow{x}=(x_{1},x_{2},x_{3},x_{4}),$ $\overrightarrow{y}%
=(y_{1},y_{2},y_{3},y_{4}),$ $\overrightarrow{z}=(z_{1},z_{2},z_{3},z_{4})$
of $\mathbb{E}^{4}$ is defined by as follows:%
\begin{eqnarray*}
\overrightarrow{x}\times \overrightarrow{y}\times \overrightarrow{z}
&=&(x_{2}y_{3}z_{4}-x_{2}y_{4}z_{3}-x_{3}y_{2}z_{4}+x_{3}y_{4}z_{2}+x_{4}y_{2}z_{3}-x_{4}y_{3}z_{2},
\\
&&-x_{1}y_{3}z_{4}+x_{1}y_{4}z_{3}+x_{3}y_{1}z_{4}-x_{3}z_{1}y_{4}-y_{1}x_{4}z_{3}+x_{4}y_{3}z_{1},
\\
&&+x_{1}y_{2}z_{4}-x_{1}y_{4}z_{2}-x_{2}y_{1}z_{4}+x_{2}z_{1}y_{4}+y_{1}x_{4}z_{2}-x_{4}y_{2}z_{1},
\\
&&-x_{1}y_{2}z_{3}+x_{1}y_{3}z_{2}+x_{2}y_{1}z_{3}-x_{2}y_{3}z_{1}-x_{3}y_{1}z_{2}+x_{3}y_{2}z_{1}).
\end{eqnarray*}%
For a hypface $\mathbf{x}$ in 4-space, we have 
\begin{equation}
I\mathbf{=}\left( 
\begin{array}{ccc}
E & F & A \\ 
F & G & B \\ 
A & B & C%
\end{array}%
\right) ,\text{ }II\mathbf{=}\left( 
\begin{array}{ccc}
L & M & P \\ 
M & N & T \\ 
P & T & V%
\end{array}%
\right) ,\text{ }III\mathbf{=}\left( 
\begin{array}{ccc}
X & Y & O \\ 
Y & Z & S \\ 
O & S & U%
\end{array}%
\right) ,  \label{I II III}
\end{equation}%
$\allowbreak \allowbreak $and%
\begin{eqnarray*}
\det I &=&(EG-F^{2})C-EB^{2}+2FAB-GA^{2}, \\
\det II &=&\left( LN-M^{2}\right) V-LT^{2}+2MPT-NP^{2}, \\
\det III &=&\left( XZ-Y^{2}\right) U-ZO^{2}+2OSY-XS^{2},
\end{eqnarray*}%
where $E=\left\langle \mathbf{x}_{u},\mathbf{x}_{u}\right\rangle ,$ $%
F=\left\langle \mathbf{x}_{u},\mathbf{x}_{v}\right\rangle ,$ $G=\left\langle 
\mathbf{x}_{v},\mathbf{x}_{v}\right\rangle ,$ $A=\left\langle \mathbf{x}_{u},%
\mathbf{x}_{w}\right\rangle ,$ $B=\left\langle \mathbf{x}_{v},\mathbf{x}%
_{w}\right\rangle ,$ $C=\left\langle \mathbf{x}_{w},\mathbf{x}%
_{w}\right\rangle ,$ $L=\left\langle \mathbf{x}_{uu},\mathbf{G}\right\rangle
,$ $M=\left\langle \mathbf{x}_{uv},\mathbf{G}\right\rangle ,$ $%
N=\left\langle \mathbf{x}_{vv},\mathbf{G}\right\rangle \mathbf{,}$ $%
P=\left\langle \mathbf{x}_{uw},\mathbf{G}\right\rangle ,$ $T=\left\langle 
\mathbf{x}_{vw},\mathbf{G}\right\rangle ,$ $V=\left\langle \mathbf{x}_{ww},%
\mathbf{G}\right\rangle ,$ $X=\left\langle \mathbf{G}_{u},\mathbf{G}%
_{u}\right\rangle ,$ $Y=\left\langle \mathbf{G}_{u},\mathbf{G}%
_{v}\right\rangle ,$ $Z=\left\langle \mathbf{G}_{v},\mathbf{G}%
_{v}\right\rangle ,$ $O=\left\langle \mathbf{G}_{u},\mathbf{G}%
_{w}\right\rangle ,$ $S=\left\langle \mathbf{G}_{v},\mathbf{G}%
_{w}\right\rangle ,$ $U=\left\langle \mathbf{G}_{w},\mathbf{G}%
_{w}\right\rangle .$ Here, 
\begin{equation}
\mathbf{G}=\frac{\mathbf{x}_{u}\times \mathbf{x}_{v}\times \mathbf{x}_{w}}{%
\left\Vert \mathbf{x}_{u}\times \mathbf{x}_{v}\times \mathbf{x}%
_{w}\right\Vert }  \label{Gauss map}
\end{equation}%
is unit normal (i.e. the Gauss map) of hypface $\mathbf{x}$. On the other
hand, $I^{-1}$\textperiodcentered $II$ gives shape operator matrix $\mathbf{S%
}$ of hypface $\mathbf{x}$ in 4-space. See \cite{GulerMagidYayli,
GulerHacisalihogluKim, GulerTurgay} for details.

\section{$i$-th Curvatures and the Fourth Fundamental Form}

To compute the $i$-th mean curvature formula $\mathfrak{C}_{i}$, where $%
i=0,..,3,$ we use characteristic polynomial $P_{\mathbf{S}}(\lambda
)=a\lambda ^{3}+b\lambda ^{2}+c\lambda +d=0$:%
\begin{equation*}
P_{\mathbf{S}}(\lambda )=\det (\mathbf{S}-\lambda I_{3})=0.
\end{equation*}%
Then, get $\mathfrak{C}_{0}=1$ (by definition), $\binom{3}{1}\mathfrak{C}%
_{1}=\binom{3}{1}H=-\frac{b}{a},$ $\binom{3}{2}\mathfrak{C}_{2}=\frac{c}{a},$
$\binom{3}{3}\mathfrak{C}_{3}=K=-\frac{d}{a}.$

Therefore, we reveal $i$-th curvature folmulas depends on the coefficients
of $I$ and $II$ fundamental forms in 4-space (It also can write depends on
the coefficients of $II$ and $III,$ or $III$ and $IV)$:

\textbf{Theorem 1.} \textit{Any hypface }$\mathbf{x}$ \textit{in }$\mathbb{E}%
^{4}$\textit{\ has following curvature formulas, }$\mathfrak{C}_{0}=1$ 
\textit{(by definition),} 
\begin{eqnarray}
\mathfrak{C}_{1} &=&\frac{\left\{ 
\begin{array}{c}
(EN+GL-2FM)C+(EG-F^{2})V-LB^{2}-NA^{2} \\ 
-2(APG-BPF-ATF+BTE-ABM)%
\end{array}%
\right\} }{3\left[ (EG-F^{2})C-EB^{2}+2FAB-GA^{2}\right] },  \label{H1} \\
\mathfrak{C}_{2} &=&\frac{\left\{ 
\begin{array}{c}
\left( EN+GL-2FM\right) V+\left( LN-M^{2}\right) C-ET^{2}-GP^{2} \\ 
-2\left( APN-BPM-ATM+BTL-PTF\right)%
\end{array}%
\right\} }{3\left[ (EG-F^{2})C-EB^{2}+2FAB-GA^{2}\right] },  \label{H2} \\
\mathfrak{C}_{3} &=&\frac{\left( LN-M^{2}\right) V-LT^{2}+2MPT-NP^{2}}{%
(EG-F^{2})C-EB^{2}+2FAB-GA^{2}}.  \label{H3}
\end{eqnarray}%
Proof. Solving $\det (\mathbf{S}-\lambda I_{3})=0$ with some algebraic
computations, we obtain coefficients of polynomial $P_{\mathbf{S}}(\lambda
). $

\textbf{Corollary 1.} For \textit{any hypface }$\mathbf{x}$ \textit{in }$%
\mathbb{E}^{4}$,\textit{\ the fourth fundamental form is related by}%
\begin{equation}
\mathfrak{C}_{0}IV-3\mathfrak{C}_{1}III+3\mathfrak{C}_{2}II-\mathfrak{C}%
_{3}I=0.  \label{Fourth}
\end{equation}

Proof. Taking $n=3$ in $\left( \ref{Fundamentalforms}\right) $, it is clear.

\textbf{Definition 1. }\textit{In 4-space, for} \textit{any hypface }$%
\mathbf{x}$ \textit{with its shape operator} $\mathbf{S}$\textit{\ and the
first fundamental form }$\left( g_{ij}\right) =I$\textit{,} \textit{%
following relations holds:\ }

$\left( a\right) $ \textit{the second fundamental form }$\left(
h_{ij}\right) =II$ \textit{is given by }$II=I$\textperiodcentered $\mathbf{S}
$\textit{,}

$\left( b\right) $ \textit{the third fundamental form }$\left( e_{ij}\right)
=III$ \textit{is given by }$III=II$\textperiodcentered $\mathbf{S}$\textit{,}

$\left( c\right) $ \textit{the fourth fundamental form }$\left(
f_{ij}\right) =IV$ \textit{is given by }$IV=III$\textperiodcentered $\mathbf{%
S}$\textit{.}

\textbf{Corollary 2. }For \textit{any hypface }$\mathbf{x}$ \textit{in }$%
\mathbb{E}^{4}$\textit{, shape operator matrix has following relation}%
\begin{equation*}
I\left( \mathbf{S}^{3}-3\mathfrak{C}_{1}\mathbf{S}^{2}+3\mathfrak{C}_{2}%
\mathbf{S}-\mathfrak{C}_{3}\right) =0.
\end{equation*}

Proof. Considering Definition 1 and Corollary 2, we write $IV=III$%
\textperiodcentered $\mathbf{S=}II$\textperiodcentered $\mathbf{S}^{2}=I$%
\textperiodcentered $\mathbf{S}^{3}.$ Then, it is clear.

\textbf{Corollary 3. }\textit{In }$\mathbb{E}^{4}$\textit{, the fundamental
forms of} \textit{any hypface }$\mathbf{x}$\textit{\ is related by}%
\begin{equation*}
\mathfrak{C}_{3}=\frac{\det III}{\det II}=\frac{\det IV}{\det III}.
\end{equation*}

Proof. Since $\mathfrak{C}_{3}=K=\frac{\det II}{\det I},$ and considering
Definition 1, it can be seen, easily.

\textbf{Corollary 4. }For \textit{any hypface }$\mathbf{x}$ \textit{in }$%
\mathbb{E}^{4}$\textit{,\ the fourth fundamental form is given by}%
\begin{equation*}
IV=3\mathfrak{C}_{1}III-3\mathfrak{C}_{2}II+\mathfrak{C}_{3}I\mathbf{.}
\end{equation*}

Proof. From $\left( \ref{Fourth}\right) $\textbf{,} we get following
symmetric fourth fundamental form matrix $\left( f_{ij}\right) $ of hypface%
\textit{\ }$\mathbf{x}$ in\textit{\ }$\mathbb{E}^{4}$:%
\begin{equation*}
IV\mathbf{=}\left( 
\begin{array}{ccc}
\zeta & \eta & \delta \\ 
\eta & \phi & \sigma \\ 
\delta & \sigma & \xi%
\end{array}%
\right) .
\end{equation*}
Here, coefficients of $IV$ are as follows%
\begin{equation*}
\zeta =-\frac{1}{\det I}\left\{ 
\begin{array}{c}
CL^{2}N-CLM^{2}-2BL^{2}T-GLP^{2}+B^{2}LX \\ 
+A^{2}NX+GL^{2}V+F^{2}VX+NP^{2}E+M^{2}VE \\ 
-CNXE+2BTXE-2MPTE-GVXE-2ABMX \\ 
-2ALNP+2BLMP+2ALMT+2CFMX-CGLX \\ 
+2AGPX-2BFPX-2AFTX-2FLMV+2FLPT%
\end{array}%
\right\} ,
\end{equation*}%
\begin{equation*}
\eta =\frac{1}{\det I}\left\{ 
\begin{array}{c}
CM^{3}-2BM^{2}P-2AM^{2}T-FNP^{2}+GMP^{2} \\ 
-FLT^{2}-B^{2}LY-A^{2}NY+FM^{2}V-F^{2}VY+MT^{2}E \\ 
+CNYE-2BTYE-MNVE+GVYE+2ABMY \\ 
-CLMN+2AMNP+2BLMT-2CFMY+CGLY \\ 
-2AGPY+2BFPY+2AFTY+FLNV-GLMV%
\end{array}%
\right\} ,
\end{equation*}%
\begin{equation*}
\delta =\frac{1}{\det I}\left\{ 
\begin{array}{c}
GP^{3}-B^{2}LO-A^{2}NO+ANP^{2}-2BMP^{2}+CM^{2}P \\ 
-ALT^{2}-AM^{2}V-2FP^{2}T-F^{2}OV+PT^{2}E+CNOE \\ 
-2BOTE+GOVE-NPVE+2ABMO-2CFMO \\ 
+CGLO-2AGOP+2BFOP+2AFOT-CLNP \\ 
+ALNV+2BLPT+2FMPV-GLPV%
\end{array}%
\right\} ,
\end{equation*}%
\begin{equation*}
\phi =-\frac{1}{\det I}\left\{ 
\begin{array}{c}
CLN^{2}-CM^{2}N-2AN^{2}P+GLT^{2}+B^{2}LZ \\ 
+A^{2}NZ+GM^{2}V+F^{2}VZ-NT^{2}E+N^{2}VE \\ 
-CNZE+2BTZE-GVZE-2ABMZ+2BMNP \\ 
+2AMNT-2BLNT+2CFMZ-CGLZ+2AGPZ \\ 
-2BFPZ-2AFTZ-2FMNV+2FNPT-2GMPT%
\end{array}%
\right\} ,
\end{equation*}%
\begin{equation*}
\sigma =\frac{1}{\det I}\left\{ 
\begin{array}{c}
T^{3}E-BNP^{2}-B^{2}LR-A^{2}NR-2AMT^{2}+BLT^{2} \\ 
+CM^{2}T-BM^{2}V-2FPT^{2}+GP^{2}T-F^{2}RV \\ 
+CNRE-2BRTE+GRVE-NTVE+2ABMR \\ 
-2CFMR+CGLR-2AGPR+2BFPR+2AFRT \\ 
-CLNT+BLNV+2ANPT+2FMTV-GLTV%
\end{array}%
\right\} ,
\end{equation*}%
\begin{equation*}
\xi =-\frac{1}{\det I}\left\{ 
\begin{array}{c}
CNP^{2}+B^{2}LS+A^{2}NS+CLT^{2}-2FMV^{2}+GLV^{2} \\ 
-GP^{2}V+F^{2}SV+NV^{2}E-T^{2}VE-CNSE \\ 
+2BSTE-GSVE-2ABMS+2CFMS \\ 
-CGLS+2AGPS-2BFPS-2AFST-2CMPT \\ 
-2ANPV+2BMPV+2AMTV-2BLTV+2FPTV%
\end{array}%
\right\} ,
\end{equation*}%
and $\det I=(EG-F^{2})C-EB^{2}+2FAB-GA^{2}$.

\subsection{$i$-th Curvatures and Fundamental Forms of Rotational
Hypersurface}

We consider the $i$-th curvatures of the rotational hypersurface $\left( \ref%
{1}\right) $, that is%
\begin{equation}
\mathbf{x}(u,v,w)=\left( f\left( u\right) \cos v\cos w,f\left( u\right) \sin
v\cos w,f\left( u\right) \sin w,\varphi (u)\right) ,  \label{6}
\end{equation}%
where $u\in \mathbb{R-\{}0\mathbb{\}}$ and $0\leq v,w\leq 2\pi .$Then, we
obtain $i$-th curvatures of $\left( \ref{6}\right) .$

Using the first differentials of rot-hypface $\left( \ref{6}\right) $, we
get the first quantities%
\begin{equation}
I=\text{diag}\left( W,f^{2}\cos ^{2}w,f^{2}\right) ,  \label{I}
\end{equation}%
where $W=f^{\prime 2}+\varphi ^{\prime 2},$ $f=f(u),$ $f^{\prime }=\frac{df}{%
du},$ $\varphi =\varphi (u),$ $\varphi ^{\prime }=\frac{d\varphi }{du}.$ The
Gauss map of the rot-hypface is%
\begin{equation}
\mathbf{G}=\left( \frac{\varphi ^{\prime }}{W^{1/2}}\cos v\cos w,\frac{%
\varphi ^{\prime }}{W^{1/2}}\sin v\cos w,\frac{\varphi ^{\prime }}{W^{1/2}}%
\sin w,-\frac{f^{\prime }}{W^{1/2}}\right) .  \label{G}
\end{equation}%
With the second differentials and $\mathbf{G}$ of hypface $\left( \ref{6}%
\right) $, we have the second quantities%
\begin{equation}
II=\text{diag}\left( -\frac{f^{\prime }\varphi ^{\prime \prime }-f^{\prime
\prime }\varphi ^{\prime }}{W^{1/2}},-\frac{f\varphi ^{\prime }}{W^{1/2}}%
\cos ^{2}w,-\frac{f\varphi ^{\prime }}{W^{1/2}}\right) ,  \label{II}
\end{equation}%
Using the first differentials of $\left( \ref{G}\right) $, we find the third
fundamental form matrix%
\begin{equation}
III=\text{diag}\left( \frac{\left( f^{\prime }\varphi ^{\prime \prime
}-f^{\prime \prime }\varphi ^{\prime }\right) ^{2}}{W^{2}},\frac{\varphi
^{\prime 2}}{W}\cos ^{2}w,\frac{\varphi ^{\prime 2}}{W}\right) .  \label{III}
\end{equation}%
We calculate $I^{-1}.II\mathbf{,}$ then obtain shape operator matrix%
\begin{equation}
\mathbf{S}=\text{diag}\left( -\frac{f^{\prime }\varphi ^{\prime \prime
}-f^{\prime \prime }\varphi ^{\prime }}{W^{3/2}},-\frac{\varphi ^{\prime }}{%
fW^{1/2}},-\frac{\varphi ^{\prime }}{fW^{1/2}}\right) .  \label{S}
\end{equation}%
Finally, with all findings, we calculate $i$-th curvatures $\mathfrak{C}_{i}$
of rot-hypface $\left( \ref{6}\right) ,$ and them give in the following
theorem.

\textbf{Theorem 2.} \textit{Rot-hypface }$\left( \ref{6}\right) $\textit{\
has following curvatures}%
\begin{eqnarray}
\mathfrak{C}_{0} &=&1\text{ \textit{(by definition),}}  \notag \\
\mathfrak{C}_{1} &=&-\frac{2\varphi ^{\prime }W+f\left( f^{\prime }\varphi
^{\prime \prime }-f^{\prime \prime }\varphi ^{\prime }\right) }{3fW^{3/2}},
\label{H1rot} \\
\mathfrak{C}_{2} &=&\frac{\varphi ^{\prime }\left( 2f^{2}\left( f^{\prime
}\varphi ^{\prime \prime }-f^{\prime \prime }\varphi ^{\prime }\right)
-W^{^{3/2}}\right) }{3f^{3}W^{2}},  \label{H2rot} \\
\mathfrak{C}_{3} &=&-\frac{\varphi ^{\prime 2}\left( f^{\prime }\varphi
^{\prime \prime }-f^{\prime \prime }\varphi ^{\prime }\right) }{f^{2}W^{5/2}}%
.  \label{H3rot}
\end{eqnarray}%
Therefore, we have following corollaries:

\textbf{Corollary 6.} \textit{Rot-hypface }$\left( \ref{6}\right) $\textit{\
is 1-minimal iff}%
\begin{equation*}
\varphi =\mp ic_{1}^{-1/4}\text{EllipticF}\left[ i\sinh ^{-1}\left(
ic_{1}^{1/4}f\right) ,-1\right] +c_{2},
\end{equation*}%
\textit{where} $i=\left( -1\right) ^{1/2}$\textit{,} EllipticF$\left[ \phi ,m%
\right] =\int\limits_{0}^{\phi }\left( 1-m\sin ^{2}\theta \right)
^{-1/2}d\theta $ \textit{is elliptic integral,} $\phi \in \left[ -\pi /2,\pi
/2\right] $, $0\neq c_{1},c_{2}$ \textit{are constants.}

Proof. \ Solving differential eq. $2\varphi ^{\prime }W+f\left( f^{\prime
}\varphi ^{\prime \prime }-f^{\prime \prime }\varphi ^{\prime }\right) =0,$
we find solutions.\textit{\ }

\textbf{Corollary 7.} \textit{Rot-hypface }$\left( \ref{6}\right) $\textit{\
is 2-minimal iff}%
\begin{equation*}
\varphi =c_{1},\text{ }\varphi =\mp i\left( \frac{c_{1}\chi }{2\rho }+\frac{%
\log \left( \rho ^{1/2}\chi +2\rho f-c_{1}\right) }{2\rho ^{3/2}}\right)
+c_{2},
\end{equation*}%
\textit{where} $\chi =\left( 4\rho f^{2}-4c_{1}f+1\right) ^{1/2},$ $\rho
=c_{1}^{2}-1,$ $0\neq c_{1},c_{2}$ \textit{are constants}$.$

Proof. \ Solving differential eq. $2f^{2}\varphi ^{\prime }\left( f^{\prime
}\varphi ^{\prime \prime }-f^{\prime \prime }\varphi ^{\prime }\right)
-\varphi ^{\prime }W^{^{3/2}}=0,$ we have solutions.\textit{\ }

\textbf{Corollary 8.} \textit{Rot-hypface }$\left( \ref{6}\right) $\textit{\
is 3-minimal iff}%
\begin{equation*}
\varphi =c_{1},\text{ }\varphi =c_{1}f+c_{2}.
\end{equation*}

Proof. \ Solving differential eq. $\varphi ^{\prime 2}\left( f^{\prime
}\varphi ^{\prime \prime }-f^{\prime \prime }\varphi ^{\prime }\right) =0,$
we get the solutions.

Next, one can see some examples in $\mathbb{E}^{4}$.

\textbf{Example 1.} \emph{Catenoidal Hypersurface. }\textit{Taking }$f\left(
u\right) =a\cosh u$ \textit{and} $\varphi \left( u\right) =au,$ \textit{where%
} $-\infty <u<\infty $\textit{,} $0\leq v,w\leq 2\pi $\textit{,} \textit{we
get}%
\begin{equation}
\mathbf{x}(u,v,w)=(a\cosh u\cos v\cos w,a\cosh u\sin v\cos w,a\cosh u\sin
w,au).
\end{equation}%
$\mathbf{x}$ \textit{verifies }$\mathfrak{C}_{1}=-\frac{1}{3a\cosh ^{2}u},$ $%
\mathfrak{C}_{2}=-\frac{1}{3a^{2}\cosh ^{4}u},$ $\mathfrak{C}_{3}=\frac{1}{%
a^{3}\cosh ^{6}u}$.

\textbf{Example 2.} \emph{Hypersphere}\textit{. Considering }$f\left(
u\right) =r\cos u$ \textit{and} $\varphi \left( u\right) =r\sin u$\textit{,
where} $r>0,$ $0<u<\pi $\textit{,} $0\leq v,w\leq 2\pi $\textit{,} \textit{%
we have}%
\begin{equation}
\mathbf{x}(u,v,w)=\left( r\cos u\cos v\cos w,r\cos u\sin v\cos w,r\cos u\sin
w,r\sin u\right) .
\end{equation}%
$\mathbf{x}$ \textit{supplies }$\mathfrak{C}_{1}=-\frac{1}{r},$ $\mathfrak{C}%
_{2}=\frac{1}{r^{2}},$ $\mathfrak{C}_{3}=-\frac{1}{r^{3}}$\textit{.}

\textbf{Example 3.} \emph{\textit{Right Spherical H}ypercylinder}\textit{.
Taking }$f\left( u\right) =r>0$ \textit{and} $\varphi \left( u\right) =u$%
\textit{, where} $0<u<\pi $\textit{,} $0\leq v,w\leq 2\pi $\textit{,} 
\textit{\ we obtain}%
\begin{equation}
\mathbf{x}(u,v,w)=\left( r\cos v\cos w,r\sin v\cos w,r\sin w,u\right) .
\end{equation}%
$\mathbf{x}$ \textit{has }$\mathfrak{C}_{1}=-\frac{2}{3r},$ $\mathfrak{C}%
_{2}=\frac{1}{3r^{2}},$ $\mathfrak{C}_{3}=0$\textit{. So, it is 3-minimal.}

Let us see some results of the fourth fundamental form of $\left( \ref{6}%
\right) .$

\textbf{Corollary 9. }\textit{The fourth fundamental form matrix} $(f_{ij})$ 
\textit{of rot-hypface\ }$\left( \ref{6}\right) $\textit{\ is as follows}%
\begin{equation}
IV=\text{diag}\left( -\frac{\left( f^{\prime }\varphi ^{\prime \prime
}-f^{\prime \prime }\varphi ^{\prime }\right) ^{3}}{W^{7/2}},-\frac{\varphi
^{\prime 3}}{fW^{3/2}}\cos ^{2}w,-\frac{\varphi ^{\prime 3}}{fW^{3/2}}%
\right) .  \label{IV}
\end{equation}

Proof. Using Corollary 4 with rot-hypface\textit{\ }$\left( \ref{6}\right) $%
, we find the fourth fundamental form matrix.

\textbf{Corollary 10. }\textit{When the curve} $\left( \ref{gamma}\right) $ 
\textit{of\ }$\left( \ref{6}\right) $ \textit{is parametrized by the arc
length, i.e.} $W=1$\textit{, then }$\left( \ref{6}\right) $ \textit{has
following relations}%
\begin{eqnarray}
\mathfrak{C}_{1}\left( 9f^{4}\mathfrak{C}_{2}+6f^{3}f^{\prime \prime
}\right) &=&2-ff^{\prime \prime }-2f^{\prime 2},  \label{c1} \\
3f^{3}\mathfrak{C}_{1}\mathfrak{C}_{3} &=&-2f^{\prime \prime }+2f^{\prime
2}f^{\prime \prime }+ff^{\prime \prime 2},  \label{c2} \\
-3ff^{\prime \prime }\mathfrak{C}_{2}-\mathfrak{C}_{3} &=&2f^{\prime \prime
2}  \label{c3}
\end{eqnarray}

Proof. The curvatures $\left( \ref{H1rot}\right) ,\left( \ref{H2rot}\right) $
and $\left( \ref{H3rot}\right) $ of the rot-hypface\textit{\ }$\left( \ref{6}%
\right) $\textit{\ }reduces to \textit{\ }%
\begin{equation}
\mathfrak{C}_{0}=1,\text{ }\mathfrak{C}_{1}=\frac{-2+2f^{\prime
2}+ff^{\prime \prime }}{3f\varphi ^{\prime }},\text{ }\mathfrak{C}_{2}=\frac{%
-2f^{2}f^{\prime \prime }-\varphi ^{\prime }}{3f^{3}},\text{ }\mathfrak{C}%
_{3}=\frac{\varphi ^{\prime }f^{\prime \prime }}{f^{2}}.
\end{equation}

\textbf{Corollary 11. }\textit{When} $f=u\neq 0$ \textit{in the previous
corollary, then }$\left( \ref{6}\right) $\textit{\ has following results,
respectively,}%
\begin{eqnarray}
&&1\text{\textit{-minimal or }}2\text{\textit{-minimal, if} }\left( \ref{c1}%
\right) \text{ \textit{holds,}} \\
&&1\text{\textit{-minimal} \textit{or} }3\text{\textit{-minimal, if} }\left( %
\ref{c2}\right) \text{ \textit{holds,}} \\
&&3\text{\textit{-minimal, if} }\left( \ref{c3}\right) \text{ \textit{holds.}%
}
\end{eqnarray}

Proof. Taking $f=u,$ it is clear.

\section{The Fourth Laplace-Beltrami Operator of a Hypersurface}

\label{Lap}The fourth Laplace-Beltrami operator of a smooth function $\phi
=\phi (x^{1},x^{2},x^{3})_{\mid _{\mathbf{D}}}$ $(\mathbf{D}\subset {\mathbb{%
R}}^{3})$ of class $C^{3}$ with respect to the fourth fundamental form of
hypface $\mathbf{x}$ is the operator $\Delta ^{IV}$, defined by%
\begin{equation}
\Delta ^{IV}\phi =\frac{1}{\mathfrak{f}^{1/2}}\dsum_{i,j=1}^{3}\frac{%
\partial }{\partial x^{i}}\left( \mathfrak{f}^{1/2}f^{ij}\frac{\partial \phi 
}{\partial x^{j}}\right) .  \label{11}
\end{equation}%
where $\left( f^{ij}\right) =\left( f_{ij}\right) ^{-1}$ and%
\begin{eqnarray*}
\mathfrak{f} &=&\det \left( f_{ij}\right) \\
&=&f_{11}f_{22}f_{33}-f_{11}f_{23}f_{32}-f_{12}f_{21}f_{33}+f_{12}f_{31}f_{23}+f_{21}f_{13}f_{32}-f_{13}f_{22}f_{31}.
\end{eqnarray*}

\subsection{Rotational Hypersurfaces Satisfying $\Delta ^{IV}\mathbf{x=Ax}$}

We now consider rot-hypface $\left( \ref{6}\right) $, and $\left( \ref{IV}%
\right) $ with $\left( \ref{11}\right) $, then have following theorem:

\textbf{Theorem 3.} \textit{The fourth Laplace-Beltrami operator of
rot-hypface} $\left( \ref{6}\right) $ \textit{is related by }$\Delta ^{IV}%
\mathbf{x}=\mathbf{Ax}$\textit{, where }$\mathbf{A}=$diag$\left( \Omega
_{1},\Omega _{2},\Omega _{3},\Phi \right) $\textit{, }%
\begin{eqnarray}
\frac{f^{3}W^{13/4}}{\varphi ^{\prime 3}\psi ^{3/2}}\left\{ f^{\prime }\frac{%
\partial }{\partial u}\left( \frac{\varphi ^{\prime 3}W^{1/4}}{f\psi ^{3/2}}%
\right) +\frac{f^{\prime \prime }\varphi ^{\prime 3}W^{1/4}}{f\psi ^{3/2}}-2%
\frac{f\psi ^{3/2}}{W^{7/4}}\right\} &=&\Omega _{i}f\text{\textit{,}}
\label{om} \\
\frac{fW^{13/4}}{\varphi ^{\prime 3}\psi ^{3/2}}\left( \varphi ^{\prime }%
\frac{\partial }{\partial u}\left( \frac{\varphi ^{\prime 3}W^{1/4}}{f\psi
^{3/2}}\right) +\varphi ^{\prime \prime }\right) &=&\Phi \varphi \text{%
\textit{,}}  \label{ph}
\end{eqnarray}%
\textit{and} $W=f^{\prime 2}+\varphi ^{\prime 2},$ $\psi =f^{\prime }\varphi
^{\prime \prime }-f^{\prime \prime }\varphi ^{\prime }$\textit{,} $i=1,2,3$%
\textit{.}

Proof. By using $\left( \ref{6}\right) ,\left( \ref{IV}\right) $ with $%
\left( \ref{11}\right) $, we compute%
\begin{equation*}
\Delta ^{IV}\mathbf{x}=\frac{1}{\left\vert \mathfrak{f}\right\vert ^{1/2}}%
\left\{ \frac{\partial }{\partial u}\left( \frac{f_{22}f_{33}}{\left\vert 
\mathfrak{f}\right\vert ^{1/2}}\mathbf{x}_{u}\right) -\frac{\partial }{%
\partial v}\left( -\frac{f_{11}f_{33}}{\left\vert \mathfrak{f}\right\vert
^{1/2}}\mathbf{x}_{v}\right) +\frac{\partial }{\partial w}\left( \frac{%
f_{11}f_{22}}{\left\vert \mathfrak{f}\right\vert ^{1/2}}\mathbf{x}%
_{w}\right) \right\} ,
\end{equation*}%
where $\mathfrak{f}=\det IV.$ If we assume that rot-hypface $\mathbf{x}$ is
constructed with component functions which are eigenfunctions of its
Laplacian, we will have that $\Delta ^{IV}\left( f\cos v\cos w\right)
=\Omega _{1}f\cos v\cos w$, $\Delta ^{IV}\left( f\sin v\cos w\right) =\Omega
_{2}f\sin v\cos w$, $\Delta ^{IV}\left( f\sin w\right) =\Omega _{3}f\sin w$, 
$\Delta ^{IV}\left( \varphi \right) =\Phi \varphi $. Hence,\textit{\ }$f(u)$ 
$\cos v\cos w$, $f(u)\sin v\cos w$\textit{\ }and\textit{\ }$f(u)\sin w$ are
eigenfunctions of $\Delta ^{IV}$ for $\Omega _{1},\Omega _{2},\Omega _{3},$
respectively, iff $f(u)$ supplies $\left( \ref{om}\right) .$ So, $\Omega
_{1}=\Omega _{2}=\Omega _{3}$ $(=\Omega $ for short). Additionally, $\varphi
\left( u\right) $ is an eigenfunction with eigenvalue $\Phi $ of $\Delta
^{IV}$ iff $\left( \ref{ph}\right) $ holds.

\textbf{Acknowledgement} \textit{The authors would like to thank the
referees for their valuable suggestions and critical remarks for improving
this paper.}

\end{document}